\documentclass[reqno,12pt]{amsart}

\usepackage{amsmath}
\usepackage{amssymb}
\usepackage{amsthm}
\usepackage[english]{babel}
\usepackage{graphics, xcolor}
\newtheorem{theorem}{Theorem}

\newtheorem{lemma}{Lemma}

\newcommand{\beqa}{\begin{eqnarray}}
\newcommand{\beqan}{\begin{eqnarray*}}
\newcommand{\eeqa}{\end{eqnarray}}
\newcommand{\eeqan}{\end{eqnarray*}}
\def\beq#1\eeq{\begin{equation}#1\end{equation}}

 \def\na{\,\, {\raise.4pt\hbox{$\shortmid$}}{\hskip-2.0pt\to}\, \, }

\def\={\overset{ \text{\rm def} }=}

\newcommand{\tc}{}\newcommand{\bl}{}

\begin{document}

\title[Concentration Functions]
{Arak Inequalities for Concentration Functions\\
and  the Littlewood--Offord Problem:\\
a shortened version}

\author[Yu.S. Eliseeva]{Yulia S. Eliseeva}
\author[F. ~G\"otze]{Friedrich G\"otze}
\author[A.Yu. Zaitsev]{Andrei Yu. Zaitsev}

\email{pochta106@yandex.ru}
\address{St.~Petersburg State University\bigskip}
\email{goetze@math.uni-bielefeld.de}
\address{Fakult\"at f\"ur Mathematik,\newline\indent
Universit\"at Bielefeld, Postfach 100131,\newline\indent D-33501 Bielefeld,
Germany\bigskip}
\email{zaitsev@pdmi.ras.ru}
\address{St.~Petersburg Department of Steklov Mathematical Institute
\newline\indent
Fontanka 27, St.~Petersburg 191023, Russia\newline\indent
and St.~Petersburg State University}

\begin{abstract}
Let $X,X_1,\ldots,X_n$ be independent identically distributed random variables.
In this paper we study the behavior of concentration functions of
 weighted sums $\sum_{k=1}^{n} X_k a_k$ with respect to the
arithmetic structure of coefficients~$a_k$ in the context of the
Littlewood--Offord problem. Concentration results of this type
received renewed interest in connection with distributions of
singular values of random matrices. Recently, Tao and Vu proposed
an Inverse Principle in the Littlewood--Offord problem. We discuss
the relations between the Inverse Principle of Tao and Vu as well
as that of Nguyen and Vu and \tc{a similar principle formulated
{\it for sums of  arbitrary independent random variables} in the
work  of Arak from the 1980's.}
\end{abstract}

\footnotetext[1]{The first author is supported by Laboratory of
 Chebyshev in St. Petersburg State University (grant of the Government of Russian Federation 11.G34.31.0026)
 and
by grant of St. Petersburg State University
6.38.672.2013.}\footnotetext[2]{The first and the third authors
are supported by grants RFBR 13-01-00256 and NSh-2504.2014.1.}
\footnotetext[3]{The second and the third authors are supported by
the SFB 701 in Bielefeld.} \footnotetext[4]{The third author was
supported by the Program of Fundamental Researches of Russian
Academy of Sciences "Modern Problems of Fundamental Mathematics".}



\maketitle

This paper is a shortened and edited version of the preprint
\cite{EGZ2}. Here we present the results without proofs.

At the beginning of 1980's, Arak \cite{Arak0, Arak} has published
new bounds for the concentration functions of sums of independent
random variables. These bounds were formulated in terms of the
arithmetic structure of supports of distributions of summands.
Using these results, he has obtained the final solution of an old
problem posed by Kolmogorov~\cite{K}. In this paper, we apply
Arak's results to
 the Littlewood--Offord
problem which was intensively investigated in the last years. We
compare the consequences of Arak's results with recent results of
Nguyen, Tao and Vu \cite{Nguyen and Vu}, \cite{Nguyen and Vu13}
and~\cite{Tao and Vu}.

Let $X,X_1,\ldots,X_n$ be independent identically distributed
(i.i.d.) $\mathbf{R}$-valued random variables. The
 concentration function of a $\mathbf{R}^d$-dimensional random
vector $Y$ with distribution $F=\mathcal L(Y)$ is defined by the
equality
\begin{equation}
Q(F,\lambda)=\sup_{x\in\mathbf{R}^d}\mathbf{P}(Y\in x+ \lambda B),
\quad \lambda\geq0, \nonumber
\end{equation}
where $B=\{x\in\mathbf{R}^d:\|x\|\leq 1/2\}$ is the centered
Euclidean ball of radius 1/2.  Let $a=(a_1,\ldots,a_n)\ne 0$,
where $a_k=(a_{k1},\ldots,a_{kd})\in \mathbf{R}^d$, $k=1,\ldots,
n$. Starting with seminal papers of Littlewood and Offord
\cite{LO} and Erd\"os~\cite{Erd}, the behavior of the
concentration functions of the weighted sums
$S_a=\sum\limits_{k=1}^{n} X_k a_k$  is studied intensively. In
the sequel, let $F_a$ denote the distribution of the sum $S_a$.

In the last ten years, refined concentration results for the
weighted sums $S_a$ play an important role in the study of
singular values of random matrices (see, for instance, Nguyen and
Vu \cite{Nguyen and Vu}, Rudelson and Vershynin \cite{Rudelson and
Vershynin08, Rudelson and Vershynin}, Tao and Vu \cite{Tao and Vu,
Tao and Vu2}, Vershynin~\cite{Vershynin}, and the authors of the
present paper   \cite{Eliseeva}, \cite{Eliseeva and Zaitsev}, and
\cite{EGZ}).
 These results reflect the dependence
of the bounds on the arithmetic structure of coefficients~$a_k$
under various conditions on the vector $a\in {(\mathbf{R}^d)}^n$
and on the distribution~$\mathcal L(X)$.

Several years ago, Tao and Vu \cite{Tao and Vu} and Nguyen and Vu
\cite{Nguyen and Vu} proposed the so-called Inverse Principles in
the Littlewood--Offord problem. In the present paper, we discuss
the relations between these Inverse Principles and similar
principles formulated  by Arak (see \cite{Arak0} and \cite{Arak})
in his work form the 1980's.
\medskip

\tc{Apparently the authors of the  publications mentioned above were not
aware  of the principles introduced in
Arak \cite{Arak0} and \cite{Arak}.
Although not named as such, his Inverse Principle} is related to general bounds for concentration functions
of distributions of sums of independent one-dimensional random
variables. The results were used for the estimation of the rate of
approximation of $n$-fold convolutions of probability
distributions by infinitely divisible ones. Later, the methods
based on \tc{Arak's  Inverse Principle admitted} to prove a number of other
important results concerning infinitely divisible approximation of
convolutions of probability measures. In 1986, Arak and Zaitsev
have published a monograph \cite{Arak and Zaitsev} containing {the}
aforementioned results and a discussion of \tc{the underlying} Inverse Principle.

The text on p. 10 of \cite{Arak and Zaitsev} is an analogue of
descriptions  of the Inverse Principles in the papers of Nguyen,
Tao and Vu \cite{Nguyen and Vu}, \cite{Nguyen and Vu13}
and~\cite{Tao and Vu}. A difference being that they restrict
themselves to the classical Littlewood--Offord problem while
discussing the arithmetic structure of the coefficients
$a_1,\ldots,a_n$ only. This means that one deals with
distributions of sums of non-identically distributed random
vectors of a special type. A further difference is that, in
\cite{Nguyen and Vu} and \cite{Nguyen and Vu13}, the multivariate
case is studied as well.

Nevertheless, there are some {\it consequences} of Arak's results
which provide some Inverse Principles for the Littlewood--Offord
problem too. Some of them  have a non-empty intersection with the
results of Nguyen, Tao and Vu \cite{Nguyen and Vu, Nguyen and
Vu13, Tao and Vu, Tao and Vu3} (see Theorem~\ref{nthm8}).
Moreover, in \cite{Arak and Zaitsev}, there are some structural
results which would be apparently new in the Littlewood--Offord
problem (see Theorems~\ref{nthm4} and~\ref{thm5}) and have no
analogues in the literature. We would like to emphasize that there
are of course also some results from \cite{Nguyen and Vu, Nguyen
and Vu13, Tao and Vu, Tao and Vu3} which do not follow from the
results of Arak.
\medskip

 We denote by
$[B]_\tau$ the closed $\tau$-neighborhood of a set $B\in\mathbf
R$.  For a finite set $K$, we denote by $|K|$ the number of
elements~$x\in K$. \bl {The symbol $\times$} is used for the
direct product of sets.

 In the sequel we use the notation
$G=\mathcal{L}(X_1-X_2)$. For $\delta\ge0$, we denote
\begin{equation}p(\delta)= G\big\{\{z:|z| >
\delta\}\big\}\big\}.\end{equation} Below we will use  the
condition
\begin{equation}
G\{\{x\in\mathbf{R}:C_1<|x| < C_2\}\}\ge C_3, \label{1s8}
\end{equation}
\tc{where the values of $C_1, C_2, C_3$ will be specified in the
formulations below}.

\medskip

 For
 $\lambda\ge0$, introduce the
distribution $H^{\lambda}$, with the characteristic function
\begin{equation} \label{11}\widehat{H}^{\lambda}(t)
=\exp\Big(-\frac{\,{\lambda}\,}2\;\sum_{k=1}^{n}\big(1-\cos\left\langle
\, t,a_k\right\rangle\big)\Big),\quad t\in \mathbf{R}^d.
\end{equation} We should note that $H^{\lambda}$, is a
symmetric infinitely divisible distribution with the L\'evy
 spectral measure $M_\lambda=\frac{\,\lambda\,}4\;M^*$, where
 $M^*=\sum_{k=1}^{n}\big(E_{a_k}+E_{-a_k}\big)$, and \,$E_b$ \,is the
distribution concentrated at a point $b\in \mathbf{R}^{d}$.

Lemma~\ref{lm42} below represents an obvious bridge between the
Littlewood--Offord problem and general bounds for concentration
functions, in particular  Arak's results.

\begin{lemma}\label{lm42}
For any $\varkappa, \delta>0$,  $\tau\ge0$, we have
\begin{equation}\label{1166}
Q(F_a, \tau) \le c_1(d)\big(1+\lfloor\varkappa/\delta\rfloor
\big)^d\, Q(H^{p(\tau/\varkappa)}, \delta),
\end{equation}where $c_1(d)$ depends on $d$ only, and $\lfloor x\rfloor$ is the
largest integer~$k$ that satisfies the inequality $k< x$.
\end{lemma}

 In a recent paper of
Eliseeva and Zaitsev \cite{Eliseeva and Zaitsev2}, a more general
statement than Lemma \ref{lm42} is obtained. It gives useful
bounds if $p(\tau/\varkappa)$ is small, even if
$p(\tau/\varkappa)=0$. Letting $\delta \to0$, we see that \eqref{1166} implies that
\begin{equation}\label{11366}
Q(F_a, 0) \le c_1(d)\, Q(H^{p(0)}, 0)=c_1(d)\, H^{p(0)}\{\{0\}\},
\end{equation}
see Zaitsev \cite{Zaitsev2} for details.

 In the monograph \cite{Arak and Zaitsev}, it is also shown
 that if the concentration function of a one-dimensional
  infinitely divisible distribution is large enough, then
 the corresponding L\'evy spectral measure is concentrated
 \tc{approximately on a set with a special arithmetic structure up to a difference of}
 small measure (see Theorems~3.3 and~4.3 of
Chapter~II in \cite{Arak and Zaitsev}). Coupled with
Lemma~\ref{lm42}, these results provide  bounds in the
Littlewood--Offord problem, see Theorems~\ref{thm7}--\ref{thm5}.

 \bigskip

A set $K\subset\mathbf R^d$  is a symmetric Generalized Arithmetic
Progression (GAP) of rank $r$ if it can be expressed in the form
$$K = \big\{m_1g_1 + \cdots + m_rg_r:-L_j \le m_j \le L_j,
 \ m_j\in{\mathbf Z} \hbox{ for all }1 \le j\le  r\big\}
$$
for some $g_1,\ldots, g_r\in\mathbf R^d$; $L_1,\ldots, L_r>0$. The
numbers $g_j$ are the generators of $K$, and $\hbox{Vol}(K) =
\prod_{j=1}^{r}(2\lfloor L_j\rfloor+1)$ is the volume of $K$ (see
\cite{Nguyen and Vu}, \cite{Nguyen and Vu13} and~\cite{Tao and
Vu}).

For any positive integers $r,m\in{\mathbf N}$ we define
$\mathcal{K}_{r,m}^{(d)}$ as the collection of all symmetric GAPs
of rank $\le r$ and of volume $\le m$.

We have to mention that, instead of $\mathcal{K}_{r,m}^{(1)}$,
Arak \cite{Arak} has considered one-dimensional projections of
sets of integer points of cardinality $\le m$ contained in
symmetric convex subsets of~${\mathbf R}^r$. However, it may be
shown that each of these projections of rank $r$ and of volume $m$
may be imbedded into a one-dimensional
 symmetric GAP of rank
$r$ and of volume $\le c_2(r)\,m$ with $c_2(r)$ depending on $r$
only. Therefore, Arak's results may be easily reformulated in
terms of $\mathcal{K}_{r,m}^{(1)}$.

For any Borel measure $W$ on ${\mathbf R}$ and $\tau\ge0$ we define $\beta_{r,m}(W, \tau)$ by
\begin{equation}
\beta_{r,m}(W, \tau)=\inf_{K\in\mathcal{K}_{r,m}^{(1)}
}W\{{\mathbf R}\setminus[K]_\tau\}. \label{1s3}
\end{equation}

 \begin{theorem}\label{thm7} Let $\varkappa, \delta>0$,  $\tau\ge0$, and let $X$ be a real random variable satisfying
 condition  \eqref{1s8} with $C_1=\tau/\varkappa$, $C_2=\infty$ and $C_3=p(\tau/\varkappa)>0$.
Let $d=1$, $r,m\in{\mathbf N}$. Then
\begin{equation}
Q(F_a, \tau)\le c_3(r)\big(1+\lfloor\varkappa/\delta\rfloor
\big)\,\biggl(\frac{1}{m\sqrt{\beta_{r,m}(M, \delta)}}
+\frac{1}{(\beta_{r,m}(M, \delta))^{(r+1)/2}}\biggr), \label{1sy4}
\end{equation}
where $M=\frac{p(\tau/\varkappa)}4\;M^*$,
 $M^*=\sum_{k=1}^{n}\big(E_{a_k}+E_{-a_k}\big)$ and where $c_3(r)$ depends on $r$ only.
\end{theorem}

In order to prove Theorem~\ref{thm7}, it suffices to apply
 Lemma \ref{lm42} and Theorem~4.3 of Chapter~II in \cite{Arak and Zaitsev}.

Theorem \ref{nthm8} \tc{follows} from Theorem~\ref{thm7}. The
conditions of this theorem are weakened conditions of those used
in  the results of Nguyen, Tao and Vu \cite{Nguyen and Vu},
\cite{Nguyen and Vu13} and~\cite{Tao and Vu}.

\medskip

 \begin{theorem}\label{nthm8}
Let $X$ be a real random variable satisfying
 condition  \eqref{1s8} with $C_1=1$, $C_2=\infty$ and $C_3=p(1)>0$.
Let $d\ge1$,  $0 < \varepsilon \le 1$, $0 < \theta \le 1$, $A>0$,
$B>0$ be constants and $\tau=\tau_{n} \geq 0$ be a parameter that
may depend on $n$. Suppose that $a=(a_1,\ldots,a_n) \in
{(\mathbf{R}^d)}^n$ is a \bl{multi-subset} of\/~${\mathbf R}^d$
such that $q_j=Q(F_a^{(j)}, \tau)\ge n^{-A}$, $j=1,\ldots,d$,
where $F_a^{(j)}$ are distributions of coordinates of the
vector~$S_a$. \tc{Let $\rho_n$ denote  a non-random sequence
satisfying}
$ n^{-B}\leq\rho_n\leq1$.
 Then, for any number $n'$ between $\varepsilon n^\theta$ and $n$,
 there exists a symmetric GAP $K$ such that

 $1$. At least $n-dn'$ elements of\/ $a$ are
$\tau\rho_n$-close to $K$ in the norm $|x|= \max_j|x_j|$
 $($this means that for these elements $a_{k}$ there exist $y_{k}\in K$ such
that $\left|a_{k}-y_{k}\right|\le\tau\rho_n)$.

 $2$. $K$ has small rank $R=O(1)$, and small cardinality
 \begin{equation}
|K|\le
\prod_{j=1}^d\max\Big\{O\Big(q_j^{-1}\,\rho_n^{-1}\,(n')^{-1/2})\Big),
1\Big\}. \label{n11sp}
\end{equation}
Here we write $O(\,\cdot\,)$ if the involved constants depend on
the parameters named ``constants'' in the formulation, but not
on~$n$.
\end{theorem}

Theorem~\ref{thm7} has been proved for one-dimensional situations
and thus initially allows us
 to  prove Theorem~\ref{nthm8} for $d=1$ only.
  However, it may be shown that this one-dimensional version of
Theorem~\ref{nthm8} provides sufficiently rich  arithmetic
properties for the set $a=(a_1,\ldots,a_n) \in {(\mathbf{R}^d)}^n$
in the multivariate case as well. It suffices to apply the
one-dimensional version of Theorem~\ref{nthm8} to the
distributions $F_a^{(j)}$, $j=1,\ldots,d$.

Theorem \ref{thm7} has non-asymptotical character, it is more
general than Theorem~\ref{nthm8} and gives information about the
arithmetic structure of $a=(a_1,\ldots,a_n)$ without assumptions
like $q_j=Q(F_a^{(j)}, \tau)\ge n^{-A}$, $j=1,\ldots,d$.\medskip

\bigskip

Another one-dimensional result of Arak \cite{Arak0} allows us to
formulate another Inverse Principle type result in the
Littlewood--Offord problem.

For any $r\in{\mathbf N}$ and $u=(u_1,\ldots, u_r)\in {({\mathbf
R}^d)}^r$, $u_j\in {\mathbf R}^d$, $j=1,\ldots,r$, we introduce
the set
\begin{equation}
{K}_{1}(u)=\Big\{\sum_{j=1}^r n_j u_j:n_j\in \{-1,0,1\}, \hbox{
for }j=1,\ldots,r\Big\}.\label{1s17}
\end{equation}
 Define also a collection of sets
\begin{equation}
\mathcal{K}_{r}^{(d)}=\big\{{K}_{1}(u):u=(u_1,\ldots, u_r)\in
{({\mathbf R}^d)}^r\big\}.\label{1s15}
\end{equation}
It is easy to see that the set $K_1(u)$ is the symmetric GAP of
rank~$r$ and volume~$3^r$.

The following Theorems \ref{nthm4} and \ref{thm5} are consequences
of Theorem 3.3 of Chapter~II \cite{Arak and Zaitsev} which follows
directly from the results of Arak \cite{Arak0}.

It is interesting that, in the multivariate case, Theorems
\ref{nthm4} and \ref{thm5} are obtained by the application of
their one-dimensional versions to the distributions of coordinates
of the vector~$S_a$.

\medskip

 \begin{theorem}\label{nthm4} Let $X$ be a real random variable satisfying
 condition  \eqref{1s8} with $C_1=1$, $C_2=\infty$ and $C_3=p(1)>0$. Let $\tau_j\ge\delta_j\geq0$
 and $q_j=Q(F_a^{(j)}, \tau_j)$, $j=1,\ldots,d$.
Then there exist an absolute constant~$c$, numbers
$r_1,\ldots,r_d\in\mathbf N$ and vectors
$u^{(j)}=\big(u_1^{(j)},\ldots, u_{r_j}^{(j)}\big)\in {{\mathbf
R}^{r_j}}$, $j=1,\ldots,d$, such that
\begin{equation}
R=\sum_{j=1}^{d}r_j\le c\sum_{j=1}^{d}\bigl(\left|\log
q_j\right|+\log(\tau_j/\delta_j)+1\bigr), \label{n1ss65}
\end{equation}
and
\begin{equation}
p(1) \,M^*\{{\mathbf
R}^d\setminus\bl{\times}_{j=1}^{d}[K_{1}(u^{(j)})]_{\delta_j}\}\le
c\sum_{j=1}^{d} \bigl(\left|\log
q_j\right|+\log(\tau_j/\delta_j)+1\bigr)^3, \label{n1ss68d}
\end{equation}
where $K_{1}(u^{(j)})\in\mathcal{K}_{r_j}^{(1)}$
 and  $M^*=\sum_{k=1}^{n}\big(E_{a_k}+E_{-a_k}\big)$.

Furthermore, the set ${\times}_{j=1}^{d}K_{1}(u^{(j)})$ can be
represented as $K_{1}(u)\in\mathcal{K}_{R}^{(d)}$, $u=(u_1,\ldots,
u_R)\in {({\mathbf R}^d)}^R$. Moreover, the vectors $u_s\in
{\mathbf R}^d$, $s=1,\ldots,R$, \tc{have  only one non-zero
coordinate each}. Denote
$$
 s_0=0\quad\hbox{and}\quad s_k=\sum_{j=1}^{k}r_j, \quad k=1,\ldots,d.
$$
For $s_{k-1}<s\le s_k$, the vectors $u_s$ \tc{ are  non-zero   in
the  $k$-th coordinates only and these coordinates are equal to
the sequence of  coordinates} $u_1^{(k)},\ldots, u_{r_k}^{(k)}$ of
the vectors~$u^{(k)}$.
\end{theorem}

\begin{theorem}\label{thm5} Let $X$ be a real random variable satisfying
 condition  \eqref{1s8} with $C_1=1$, $C_2=\infty$ and $C_3=p(1)>0$.
Let $A,B>0$. Let $\tau_j\ge \delta_j\geq0$, $\tau_j/\delta_j\le
n^B$ and
 $q_j=Q(F_a^{(j)}, \tau_j)\ge n^{-A}$, for  $j=1,\ldots,d$.
Then there exist an absolute constant~$c$, numbers
$r_1,\ldots,r_d\in\mathbf N$ and vectors
$u^{(j)}=\big(u_1^{(j)},\ldots, u_{r_j}^{(j)}\big)\in {{\mathbf
R}^{r_j}}$, $j=1,\ldots,d$, such that
\begin{equation}
R=\sum_{j=1}^{d}r_j\le c\, d\,\big((A+B)\,\log n+1\big),
\label{1st65}
\end{equation}
and
\begin{equation}
p(1) \,M^*\{{\mathbf
R}^d\setminus\bl{\times}_{j=1}^{d}[K_{1}(u^{(j)})]_{\delta_j}\}\le
c\, d \,\big((A+B)\,\log n+1\big)^3, \label{1st68d}
\end{equation}
where $K_{1}(u^{(j)})\in\mathcal{K}_{r_j}^{(1)}$
 and  $M^*=\sum_{k=1}^{n}\big(E_{a_k}+E_{-a_k}\big)$.
Moreover, the description of the set
$K_{1}(u)=\bl{\times}_{j=1}^{d}K_{1}(u^{(j)})$ from the end of the
formulation of Theorem~$\ref{nthm4}$ remains true.\end{theorem}

It is easy to see that, in conditions of Theorem~\ref{thm5} with
$\tau_j=\delta_j\,n^B=\tau$, $j=1,\ldots,d$, the set $K_{1}(u)$ is
a symmetric GAP of rank $R=O(\log n)$, of volume $3^R=O(n^D)$
(with a constant~$D$) and such that at least $n-O(\log^3 n)$
elements of $a=(a_1,\ldots,a_n) \in ({\mathbf{R}^d})^n$ are
$\tau/n^B$-\tc{close} to~$K_{1}(u)$. Theorem~\ref{nthm4} provide
bounds \tc{by replacing $\log n$ by $\left|\log q\right|$ without
the assumption~$q=Q(F_a, \tau)\ge n^{-A}$}. Moreover, in
\eqref{n1ss68d} and \eqref{1st68d}, the dependence of constants on
$C_3=p(1)$ is  \bl{stated explicitly}.

Notice that if $\tau_1=\cdots=\tau_d=\tau$, then $q=Q(F_a,
\tau)\le q_j$, and $\left|\log q_j\right|\le\left|\log q\right|$,
$j=1,\ldots,d$. Moreover, there exist distributions for which the
quantity $q$ may be sufficiently smaller than $\max_j q_j$.
Consider, for instance, the uniform distribution on the boundary
of the square $\big\{x\in\mathbf R^2:|x|=1\big\}$.

Moreover, if all $\tau_j$ and $\delta_j$ are equal to zero, 
we can use \eqref{11366} instead of \eqref{1166} and replace $p(1)$ by $p(0)$ 
in \eqref{n1ss68d} and \eqref{1st68d}, see Zaitsev \cite{Zaitsev2} for details.
In this case we agree that $\log(\tau_j/\delta_j)=0$.

\bigskip

Now we compare our Theorems~\ref{nthm8}, \ref{nthm4}
and~\ref{thm5} with the results discussed in a review of Nguyen
and Vu \cite{Nguyen and Vu13} (see Theorems 7.5, 9.2 and 9.3 of
\cite{Nguyen and Vu13}). These results were obtained under the
assumption $Q(F_a, \tau)\ge n^{-A}$. This implies that
$Q(F_a^{(j)}, \tau)\ge n^{-A}$, $j=1,\ldots,d$, since
$Q(F_a^{(j)}, \tau)\ge Q(F_a, \tau)$.

 A few years ago
Tao and Vu \cite{Tao and Vu}  formulated the so-called Inverse
Principle, stating that
$$\hbox{\it A set $a=(a_1,\ldots,a_n)$ with large $Q(F_a,
0)$ must have strong additive structure.}$$ Theorem~7.5 of
\cite{Nguyen and Vu13} was obtained by Tao and Vu \cite{Tao and
Vu}. This theorem is named in \cite{Nguyen and Vu13} ``Weak
Inverse Principle''.

Later, Tao and Vu \cite{Tao and Vu3} improved the result of
Theorem~7.5 of \cite{Nguyen and Vu13}.
 Nguen and Vu \cite{Nguyen and Vu}
have extended the Inverse Principle to the continuous case
proving, in particular, Theorems~9.2 and~9.3 of \cite{Nguyen and
Vu13}.

Theorem~\ref{nthm8}
 allows us to derive Theorem~7.5 of \cite{Nguyen and Vu13} and a
one-dimensional version of the first two statements of Theorem~9.3
of \cite{Nguyen and Vu13}.

Lemma \ref{lm42} is interesting only if we assume that
$p(\tau/\varkappa)>0$. This assumption is closely related to
assumption \eqref{1s8} in Theorems~9.2 and~9.3 of \cite{Nguyen and
Vu13} (with $C_2<\infty$). We can w.l.o.g. take
 in  \eqref{1s8} $C_1=1$ and $C_3=p(1)$. Moreover, in our results, $C_2=\infty$.
We think that using Lemma \ref{lm42}, one could show that $C_2$
may be taken as $C_2=\infty$
 in Theorems~9.2 and~9.3 of \cite{Nguyen and Vu13} too.
 Note, however,
 that $p(1)$ is involved in our inequalities explicitly, in contrast with Theorems~9.2 and~9.3 of \cite{Nguyen and Vu13}.

The basic inequality of Theorem~9.3 of \cite{Nguyen and Vu13}
gives the bound
\begin{equation}
|K|\le \max\big\{O(q^{-1}(n')^{-1/2}), 1\big\}, \quad\text{where }
q=Q(F_a, \tau). \label{12sp}
\end{equation}
Inequality~\eqref{12sp} and inequality~\eqref{n11sp} of
Theorem~\ref{nthm8} (with $\rho_{n}=1$) are not only of the same
form, but their contents are almost the same, at least for
$d=1$. \tc{Attentive readers 
may notice evident differences though}. In particular, the last
item of Theorem~9.3 of \cite{Nguyen and Vu13} is absent in
Theorem~\ref{nthm8}. On the other hand, in Theorem~\ref{nthm8}, we
take $C_2=\infty$.
\medskip

Sometimes, for $d>1$, inequality \eqref{n11sp} (with $\rho_{n}=1$)
may be even stronger than inequality~\eqref{12sp}. For example, if
the vector $S_a$ has independent coordinates (this may happen if
each of the vectors~$a_j$ has only one non-zero coordinate), then
 \begin{equation}
c_4 (d)\prod_{j=1}^d q_j\le Q(F_a,\tau)\le\prod_{j=1}^d q_j
 \label{11sps}
\end{equation}
with a positive $c_4(d)$. Note, however, that we could derive a
multivariate analogue of Theorem~9.3 of \cite{Nguyen and Vu13}
from its one-dimensional version arguing precisely as in the proof
of our Theorem \ref{nthm8}. Then we get inequality \eqref{n11sp}.

 Theorem~\ref{nthm8} can be considered as an analogue of both
Theorems~9.2 and~9.3 of \cite{Nguyen and Vu13}. Comparing these
theorems, we should mention that the amount of approximating
points is sometimes a little bit smaller in Theorem~9.2 of
\cite{Nguyen and Vu13}, but, in Theorem~\ref{nthm8}, $C_2=\infty$,
and we \tc{get a variety of results by choosing} various $\rho_n$,
while in Theorem~9.2 of \cite{Nguyen and Vu13}
$\rho_n=n^{-1/2}\log n$, and in Theorem 9.3 of \cite{Nguyen and
Vu13} $\rho_n=1$.
\medskip

The assertion of Theorem~\ref{thm5} implies that, in conditions of
Theorem~9.3 of \cite{Nguyen and Vu13}, there exists a symmetric
GAP $K$ of rank $R=O(\log n)$, of volume~$3^R=O(n^D)$ and such
that at least $n-O((\log n)^3)$ elements of $a=(a_1,\ldots,a_n)
\in ({\mathbf{R}^d})^n$ are $\tau/n^B$-closed to~$K$. Moreover,
Theorem~\ref{nthm4} provide bounds with replacing $\log n$ by
$\left|\log q\right|$ without assumption~$q=Q(F_a, \tau)\ge
n^{-A}$ (recall that this assumption is absent in the conditions
of Theorem~\ref{thm7} too). Comparing with the results of
\cite{Nguyen and Vu13}, we see that in Theorem~\ref{thm5} the
exceptional set has a logarithmic size (which is much better than
$O(n)$ and $O(n^\theta)$, $0<\theta\le1$, in Theorems~9.2 and~9.3
of \cite{Nguyen and Vu13}), but this is attained at the expense of
a logarithmic growth of the rank.

\end{document}